\documentclass[letterpaper,12pt,onecolumn]{IEEEtran}

\usepackage{moreverb}
\usepackage{hyperref}
\usepackage{algorithm}
\usepackage{algpseudocode}

\usepackage{subcaption}
\DeclareCaptionSubType*{algorithm}

\DeclareCaptionLabelFormat{alglabel}{Algorithm~#2}

\RequirePackage{amstext}
\RequirePackage{amsfonts}
\RequirePackage{amssymb}
\RequirePackage{amsmath} 
\RequirePackage{bm}
\usepackage{color}
\usepackage{graphicx}  
\usepackage{booktabs}
\usepackage{multirow}
\usepackage{tabularx}

\usepackage{booktabs}
\newcommand{\ra}[1]{\renewcommand{\arraystretch}{#1}}

\newcommand{\pdmat}[1]{\mathbb{S}_{++}^{#1}}
\newcommand{\psdmat}[1]{\mathbb{S}_{+}^{#1}}
\newcommand{\reals}[1]{\mathbb{R}^{#1}}
\newcommand{\nonnegreals}[1]{\reals{#1}_{+}}
\newcommand{\posreals}[1]{\reals{#1}_{++}}
\newcommand{\myref}[0]{\text{ref}}
\newcommand{\secref}[1]{Section~\ref{#1}}
\newcommand{\figref}[1]{Figure~\ref{#1}}
\newcommand{\mydef}[0]{\triangleq}

\newcommand{\dnano}[0]{\emph{dnano}}

\title{\LARGE \bf
Optimization-Based Autonomous Racing\newline of 1:43 Scale RC Cars
}

\author{Alexander Liniger, Alexander Domahidi and Manfred Morari
\thanks{ The authors are with the Automatic Control Laboratory, ETH
Zurich, 8092 Z\"urich, Switzerland; 

Emails: \tt{liniger|domahidi|morari@control.ee.ethz.ch}}
}

\begin{document}

\maketitle
\thispagestyle{empty}
\pagestyle{empty}

\begin{abstract}
This paper describes autonomous racing of RC race cars based on mathematical optimization. Using a dynamical model of the vehicle, control inputs are computed by receding horizon based controllers, where the objective is to maximize progress on the track subject to the requirement of staying on the track and avoiding opponents. Two different control formulations are presented. The first controller employs a two-level structure, consisting of a path planner and a nonlinear model predictive controller (NMPC) for tracking. The second controller combines both tasks in one nonlinear optimization problem (NLP) following the ideas of contouring control. Linear time varying models obtained by linearization are used to build local approximations of the control NLPs in the form of convex quadratic programs (QPs) at each sampling time. The resulting QPs have a typical MPC structure and can be solved in the range of milliseconds by recent structure exploiting solvers, which is key to the real-time feasibility of the overall control scheme. Obstacle avoidance is incorporated by means of a high-level corridor planner based on dynamic programming, which generates convex constraints for the controllers according to the current position of opponents and the track layout. The control performance is investigated experimentally using 1:43 scale RC race cars, driven at speeds of more than 3 m/s and in operating regions with saturated rear tire forces (drifting). The algorithms run at 50 Hz sampling rate on embedded computing platforms, demonstrating the real-time feasibility and high performance of optimization-based approaches for autonomous racing.
\end{abstract}

\maketitle

\section{Introduction}
Autonomous car racing is a challenging task for automatic control systems due to the need for handling the vehicle close to its stability limits and in highly nonlinear operating regimes. In addition, dynamically changing racing situations require advanced path planning mechanisms with obstacle avoidance executed in real-time. Fast dynamics constrain the sampling time to be in the range of a few tens of milliseconds at most, which severely limits the admissible computational complexity of the algorithms. This situation is even more challenging if the autonomous algorithms shall be executed on simple, low-power embedded computing platforms.

In this paper, we investigate optimization-based control strategies for the task of racing an autonomous vehicle around a given track. We focus on methods that can be implemented to run in real-time on embedded control platforms, and present experimental results using 1:43 scale Kyosho \dnano~RC race cars that achieve top speeds of more than 3\,m/s, which corresponds to an upscaled speed of about 465 km/h. For high performance, the proposed controllers operate the car at its friction limits, far beyond the linear region of what is typically used in other autonomous driving systems. This challenging task is generally mastered only by expert drivers with lots of training. In contrast, our approach requires merely a map of the track and a dynamical model of the car; in particular, we use a bicycle model with nonlinear tire forces, neglecting load transfer and coupling of lateral and longitudinal slip.

Two model-based control schemes are presented in this paper. First, we describe a hierarchical two-level control scheme consisting of a model-based path planner generating feasible trajectories for an underlying nonlinear model predictive control (NMPC) trajectory tracking controller. Our second approach combines both path planning and path following into one formulation, resulting in one NMPC controller that is based on a particular formulation known from contouring control~\cite{Faulwasser2009,Lam2010}. The objective of both approaches is to maximize the progress on the track, measured by a projection of the vehicle's position onto the center line of the track. Linear time varying (LTV) models obtained by linearization are employed to construct a tractable convex optimization problem to be solved at each sampling time. Efficient interior point solvers for embedded systems, generated by FORCES~\cite{Domahidi2012,FORCEScodegen}, are employed to solve the resulting optimization problems, which makes the approaches presented in this paper amenable for use on embedded systems.

We furthermore demonstrate that both schemes can easily be extended to incorporate obstacle avoidance by adjusting the constraints of the resulting optimization problems according to the racing situation. With this mechanism in place, the avoidance trajectory is optimized to yield maximal overall progress, which translates into highly effective overtaking maneuvers that are automatically planned and executed.

\subsection{Related work}

Safe autonomous driving at moderate speeds has been demonstrated for example in the context of Autonomous Highway Systems (AHS) in 1997~\cite{Horowitz2000} within the Californian PATH programme, or in contests such as the DARPA Grand Challenge in 2005~\cite{Buehler2007} and the Urban Challenge in 2007~\cite{Buehler2009} by numerous research groups. In these projects, the main goal was do develop autonomous driving systems for public traffic situations. Autonomous racing has received less attention, but impressive advances have been made recently also in this discipline. For instance, a fully autonomous Audi TTS has been reported to race at high speeds in~\cite{stanford} using a trajectory tracking controller, which follows a pre-computed trajectory. However, to the best knowledge of the authors, fully automatic racing including obstacle avoidance in dynamic racing situations has not been demonstrated yet in practice.

From a control perspective, the fundamental building blocks for automatic racing can be grouped into three categories: drift control, reference tracking, and path planning. Research associated with the first group focuses on gaining a better understanding of the behavior of the car near the friction limit for designing safety control systems that try to prevent loss-of-control accidents. By analyzing the nonlinear car model, it can be shown that steady state motions corresponding to drifting can be generated, see e.g.~\cite{Velenis2009} and~\cite{Voser2010}. As these approaches are not designed for trajectory tracking, the controlled car drifts in a circle, which represents a steady state drift equilibrium.

Reference tracking controllers usually are designed to operate the vehicle within the linear tire force region for maximum safety and robustness. Approaches based on NMPC, which allows one to incorporate the latter constraint in a systematic manner, deal with the nonlinearities of the model either directly by a nonlinear programming (NLP) solver~\cite{Borrelli2005}, or use LTV~\cite{Falcone2007} or piece-wise affine (PWA) approximations~\cite{Besselmann2008}, resulting in a convex quadratic program (QP) or a mixed-integer QP, respectively. Approaches without optimization in real-time include nonlinear robust control~\cite{Ackermann1995} and flatness-based control~\cite{Fuchshumer2005}, which however lack the ability to incorporate constraints. One approach that is explicitly designed for operating with saturated tire forces is~\cite{stanford}, where the center of percussion (CoP) is used to design a state feedback steering controller for reference tracking. At the CoP, the rear tire forces do not influence the dynamics by definition, which allows for a simple linear steering controller.

In order to avoid obstacles, reference tracking schemes rely on a higher-level path planner. A simple point mass model in the high-level path planner is used in~\cite{Gao2010} to limit the computational complexity. This can be problematic if the generated trajectories are infeasible with respect to the car's physical dynamics, which can lead to preventable collisions~\cite{Gray2012}. The latter reference suggests to use a library of steady-state movements (trims) with parametrizable length, which are generated from a more complex dynamical model, and to connect them by maneuvers which allow the transition from one trim to another, similar to the idea of~\cite{Frazzoli2005}. This path planning has the advantage that it generates feasible reference trajectories for the low level tracking controller, and that drifting trims can be included to increase the agility of the car. However, the resulting optimization problem is a mixed-integer program that is too complex to be solved within the typical sampling times on embedded platforms. Consequently, the approach is not well suited for real-time autonomous racing.

In order to avoid relying on the feasibility of the path planner, one-level approaches have been investigated e.g. in~\cite{Velenis2007}, where optimal trajectories and the associated open-loop control inputs for rally racing maneuvers are calculated numerically in simulation. In~\cite{Gao2010}, obstacle avoidance is incorporated into the tracking controller by using an additional cost term that penalizes coming too close to the obstacle. However, the controller is reported to perform worse than its two-level counterpart, especially if the car has to diverge from the reference trajectory. A one-level approach similar to~\cite{Gao2010} is studied in~\cite{Frasch2013} in simulation, where obstacle avoidance is incorporated into the problem formulation by using a spatial reformulation and imposing constraints on the transformed variables. While in~\cite{Velenis2007,Gao2010} the solution to the NLP is obtained by a standard nonlinear solver with prohibitively long solve times,~\cite{Frasch2013} uses real-time iterations~\cite{Diehl2002b} to achieve the low computation times needed for a real-time implementation. However, it is assumed that there is only one side where the obstacle can be overtaken, which may not be the case in practical racing situations.

An interesting alternative to optimization-based methods are sampling-based methods. For example, rapidly exploring random trees (RRTs) have been investigated in \cite{RRT1} and \cite{Jeon2013} to generate time optimal trajectories through a 180$^{\circ}$ curve. The advantage is that such algorithms tend to quickly find feasible trajectories. In \cite{Jeon2013}, the differential flatness of the system is exploited to speed up the convergence of the algorithm, generating obstacle-free trajectories even in complex situations. However, the reported computation times are not yet low enough to allow for a real-time implementation.

\subsection{Contribution}

In this paper, we describe two novel autonomous, optimization-based racing controllers that incorporate obstacle avoidance, track constraints, actuator limitation and the ability for controlled drift in a systematic and straightforward way. Real-time feasibility at 50\,Hz sampling rate is demonstrated in experimental results with fast RC cars. To the best knowledge of the authors, this is the first implementation of autonomous racing controllers with obstacle avoidance on an experimental testbed.

The fundamental idea of both approaches is to use a receding horizon controller, which maximizes progress on the track within the horizon as a performance measure. This is closely related to a time optimality criterion and allows for a systematic incorporation of obstacles and other constraints. It is particularly effective for overtaking opponents, as our controllers seek for a progress-optimal solution around the obstacles. In order to deal with track constraints and obstacle avoidance situations, we represent the feasible set for the position of the car at any time instance by a slab defined by two parallel linear inequalities. This ensures tractability of the subproblems by convex programming on one hand, and incorporation of track and obstacle constraints in a systematic manner on the other hand. To deal with non-convex situations such as overtaking an obstacle on the left or right, a high-level corridor planning algorithm supplies a set of appropriate convex constraints to the controllers. This  constraint set is the result of a shortest path problem solved by dynamic programming. 

The first, hierarchical control approach presented in this paper is in principle similar to~\cite{Gray2012}, but its complexity is low enough to be fully implemented in real-time. Our path planner also uses trims, but only \emph{one} is selected for the whole prediction horizon from a set of trajectories, which represent steady-state cornering conditions of the nonlinear model, gridded for velocity and steering angle. The path planner selects the trajectory with the largest progress that does not leave the track or hits any obstacle. Moving obstacles can inherently be dealt with, as new trims are generated every sampling time. A low level nonlinear reference tracking MPC controller tracks the selected trim, as in \cite{Gray2012}, but we additionally use obstacle avoidance constraints as described above, which turns out to be very effective in practice.

The second controller is based on a model predictive contouring control (MPCC) formulation \cite{Faulwasser2009,Lam2010}, combining path generation and path tracking into one problem. The MPCC essentially plans a progress-optimal path by taking into account the (nonlinear) projection of the vehicle's position onto the center line. The resulting controller is able to plan and to follow a path which is similar to the time-optimal path in \cite{Velenis2007} when the horizon is chosen long enough. 

Both NMPC problems are approximated locally by linearizing the continuous nonlinear system dynamics around a state trajectory to obtain an LTV model. In the hierarchical approach, we linearize around the trajectory obtained from the path planner, while in the one-level approach the linearization points are given by the shifted trajectory from the previous sampling time. We then discretize by the matrix exponential, and solve the resulting quadratic program (QP) by a tailored interior point solver, generated by FORCES~\cite{Domahidi2012,FORCEScodegen}. The computation time of the solvers scale linearly with the horizon length, hence we make use of long horizons of up to 40 time steps in the MPCC, for example. After the QP has been solved, the first control input is applied to the system, and the process is repeated at the next sampling time. This approach for solving the NMPC problems corresponds to the basic version of the real-time iterations from~\cite{Diehl2002b}.

\subsection{Outline}

The paper is organized as follows. In Section~\ref{sec:model}, the dynamical model of the car and some of its basic properties are described. In Section~\ref{sec:hrhc}, we present the hierarchical control approach with separate path planner and path tracking, while the combined approach is presented in Section~\ref{sec:mpcc}. In Section~\ref{sec:dp}, we briefly discuss our method for selecting convex constraints for the two controllers, enabling obstacle avoidance in dynamic racing situations. Section~\ref{sec:results} presents the testbed setup and experimental results for both controllers.

\section*{Notation}

\newcommand{\pd}

The following notation is used throughout the paper. The set of reals is denoted by $\reals{}$, and the set of real column vectors of length $n$ is denoted by $\reals{n}$. The set of nonnegative and strictly positive reals is denoted by $\nonnegreals{}$ and $\posreals{}$, respectively. The set of positive definite matrices of size $n$ is denoted by $\pdmat{n}$ and the set of positive semi-definite matrices of size $n$ by $\psdmat{n}$. We define $\|x\|_P^2 \triangleq x^T P x$ for a vector $x\in\reals{n}$ and some positive definite matrix $P\in \pdmat{n}$. If $a\in\reals{n}$ and $b\in\reals{m}$ are column vectors or scalars, we denote their concatenation by $(a,b)\mydef [a^T\,, b^T]^T \in\reals{n+m}$.

\section{Car Model} \label{sec:model}

\subsection{Bicycle Model}

\begin{figure}[h]
    \centering
    \includegraphics[width = 0.45\textwidth]{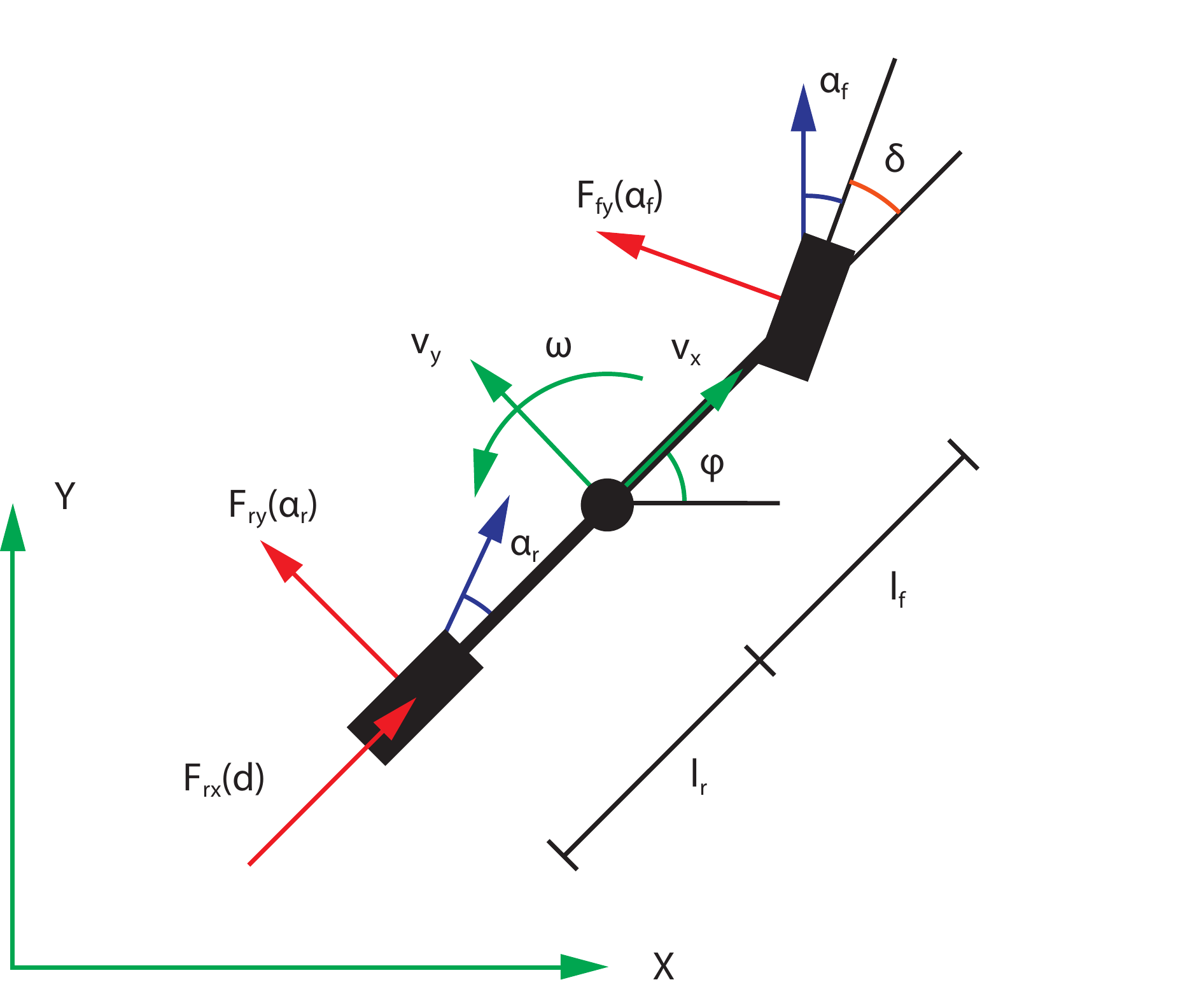}
    \caption{Schematic drawing of the car model}\label{fig:bicycle}
\end{figure}

The RC cars are modeled using a bicycle model as done in~\cite{Velenis2009,Voser2010}, where the car is modeled as one rigid body with a mass $m$ and an inertia $I_z$, and the symmetry of the car is used to reduce it to a bicycle. Only the in-plane motions are considered, i.e. the pitch and roll dynamics as well as load changes are neglected. As the used cars are rear wheel driven and do not have active brakes, the longitudinal force on the front wheel is neglected. The resulting model is shown in \figref{fig:bicycle} together with the resulting differential equations \eqref{odeModel} defining the model:

	\begin{subequations}
	\label{odeModel}
	\begin{align}
	\dot{X} &= v_x \cos(\varphi) - v_y \sin(\varphi) \,,\\
	\dot{Y} &= v_x \sin(\varphi) + v_y \cos(\varphi)\,, \\
	\dot{\varphi} &= \omega\,,\\
	\dot{v}_x &= \frac{1}{m}(F_{r,x} - F_{f,y} \sin{\delta} + m v_y \omega)\,, \\
	\dot{v}_y &= \frac{1}{m}(F_{r,y} + F_{f,y} \cos{\delta} - m v_x \omega)\,, \\
	\dot{\omega} &= \frac{1}{I_z}(F_{f,y} l_f \cos{\delta} - F_{r,y} l_r)\,.
	\end{align}
	\end{subequations}
	
The equation of motion is derived around the center of gravity (CoG), where the states are the position $X$, $Y$ of the CoG in the inertial frame and the angle of the car relative to the inertial frame, $\varphi$. This is the kinematic part of the model. The kinetic part of the model is derived around a body-fixed frame centered at the CoG, where the states are the longitudinal and lateral velocity of the car, $v_x$ and $v_y$, and finally the yaw rate $\omega$. The control inputs are the PWM duty cycle $d$ of the electric drive train motor and the steering angle~$\delta$. The subscripts $x$ and $y$ indicate longitudinal and lateral forces or velocities, while $r$ and $f$ refer to the front and rear tires, respectively. Finally, $l_f$ and $l_r$ are the distance from the CoG to the front and the rear wheel. 

The tire forces $F$ model the interaction between the car and the road and are the most important part of the dynamics. As the goal is for the cars to race, the model of the tire forces has to be realistic enough to represent the car at high speeds and its handling limits. The lateral forces $F_{f,y}$ and $F_{r,y}$ are modeled using a simplified Pacejka Tire Model \cite{MF}, see~\eqref{modelForces1} and \eqref{modelForces2}, which is a good approximation to measured friction curves in practice:
\begin{subequations}
\begin{align}
\label{modelForces1}
F_{f,y} &= D_f \sin (C_f \arctan(B_f \alpha_f)) \quad  \text{where} \quad \alpha_f = -\arctan\left(\frac{\omega l_f + v_y}{v_x}\right) + \delta\,,\\
\label{modelForces2}
F_{r,y} &= D_r \sin (C_r \arctan(B_r \alpha_r)) \quad  \:\text{where} \quad \alpha_r = \arctan\left(\frac{\omega l_r - v_y}{v_x}\right)\,,\\
\label{modelForces5}
F_{r,x} &= (C_{m1} - C_{m2} v_x) d -  C_r -  C_d v_x^2\,.
\end{align}
\end{subequations}
The parameters $B$, $C$ and $D$ define the exact shape of the semi-empirical curve. The longitudinal force of the rear wheel $F_{r,x}$ is modeled using a motor model for the DC electric motor as well as a friction model for the rolling resistance and the drag, cf.~\eqref{modelForces5}. Due to the size of the cars, it is currently not feasible to measure the wheel speed, which makes it impossible to use combined slip models such as~\cite{Velenis2009,Velenis2007}.

The model is identified using a combination of stationary and dynamic identification experiments, as reported in~\cite{Voser2010}. This allows for identifying the rear wheel combined slip effects into the lateral tire friction model~\eqref{modelForces2}. Thus the identified model is suitable to represent the car also at its handling limits, when the tire forces are saturated or close to saturation.

\subsection{Stationary Velocity Analysis} \label{sec:stationaryVelocities}
For a better understanding of the model, the stationary velocities of the model are investigated. A similar, more elaborate analysis is presented in~\cite{Velenis2009,Voser2010}.

The objective is to find points in the model where all accelerations are zero. This is done for different constant forward velocities $\bar{v}_{x}$ and constant steering angles $\bar{\delta}$, thus the problem becomes a nonlinear algebraic system of equations, with two equations ($\dot{v}_y = 0,$ $\dot{\omega}=0$) and two unknowns ($\bar{v}_{y},$ $\bar{\omega}$). By finding a solution to the algebraic equations for different steering angles $\bar{\delta}$ and different initial conditions, the stationary velocities for one forward velocity ($\bar{v}_{y},$ $\bar{\omega}$) can be calculated.

\begin{figure}
\centering
  \includegraphics[width = 0.9\textwidth]{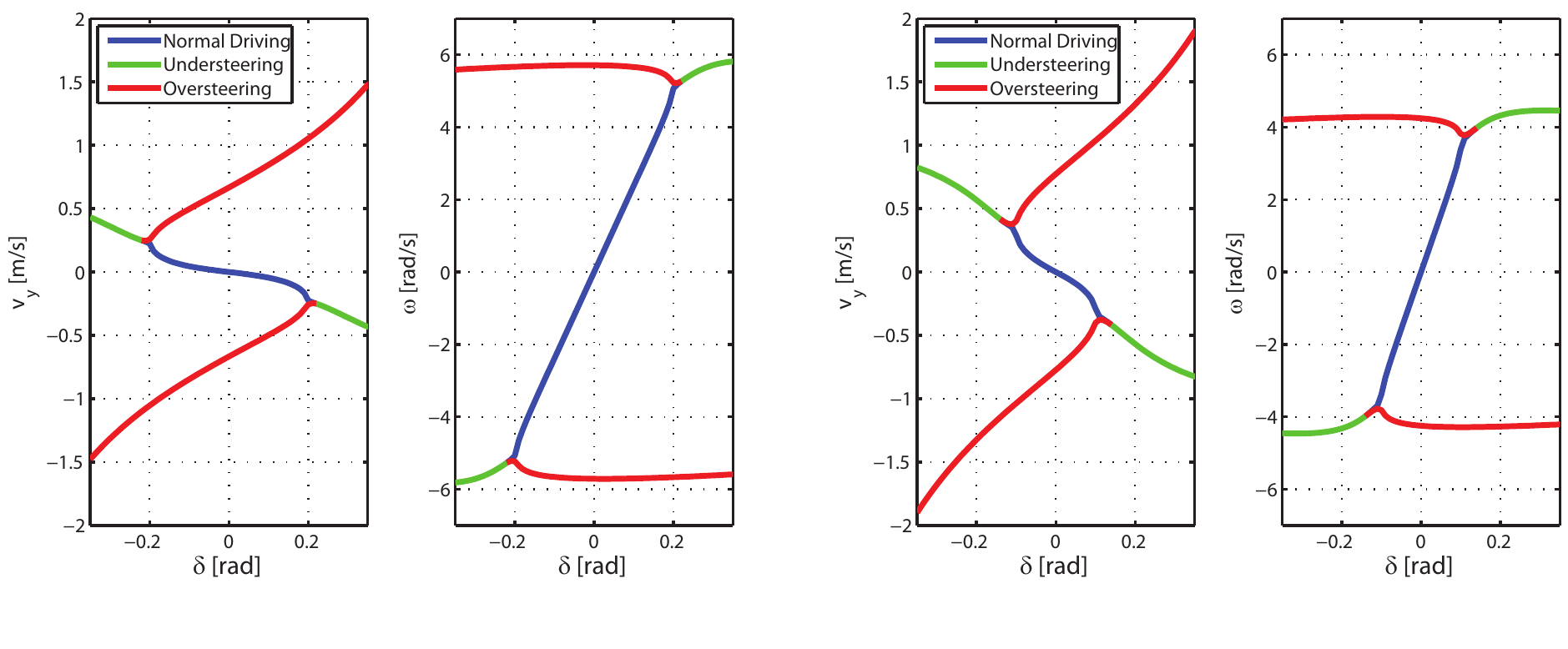}
  \caption{Stationary velocities for $\bar{v}_x=1.5$\,m/s (plots on the left) and $\bar{v}_x=2$\,m/s (plots on the right), parameterized by the steering angle $\bar{\delta}$.}
  \label{statVelo}
\end{figure}

The resulting stationary velocities are shown for $\bar{v}_x=1.5$\,m/s and $\bar{v}_x=2$\,m/s in Figure \ref{statVelo}. The resulting stationary velocities can be categorized into three different regions. The blue line is the normal driving region, which is characterized by the linear dependency of the yaw rate and the steering angle, as well as small lateral velocities. The other two regions in red and green correspond to over- and understeering points in the model. In the understeering region (green curves), the saturated front tire force law prevents the car from achieving larger yaw rates, and thus the car can not take a sharper turn. In the oversteering region (red curves), the rear tire force law is fully saturated and the car drives with a high lateral velocity, which is usually called drift.

The dynamics of the model in the different regions are different, for example in the oversteering region where the car is drifting, the steering angle can be of opposite sign than the curvature the car is driving, known as counter steering, which does not occur in the case of the normal and the understeering mode.
A stability analysis of the lateral velocity model around the different stationary velocity points also shows the difference in the dynamics: While the normal driving region and the understeering region are stable, the drifting region is unstable~\cite{Voser2010}.

From Figure~\ref{statVelo} it is visible that the normal steering region is getting smaller when the car is driving faster. Thus the maximal and minimal steering angle for which the tire forces are not saturated is reduced with increasing velocities, and at the same time the maximal yaw rate decreases. 

Because all velocities are constant, the movement of the car is a uniform circular movement, for which the relationship between the yaw rate, the radius $R$ and the curvature $\kappa$ is known,
\begin{equation}
\label{yawRateToRadius}
|\omega| = \frac{v}{R} = v \kappa \quad \text{where} \quad v = \sqrt{v_x^2 + v_y^2}\,.
\end{equation}
Consequently, all stationary velocity points correspond to the car driving in a circle with a constant radius, i.e. the model can predict drifting along a circle.

\section{Autonomous Racing Control} \label{sec:control}
 
In this main part of the paper, we present two optimization-based formulations for the task of autonomous racing for RC cars which are based on the model described in \secref{sec:model}. First, the hierarchical two-level approach is presented in \secref{sec:hrhc}, followed by the one-level scheme in \secref{sec:mpcc}. We then briefly discuss a high-level method for obstacle avoidance that is the same for both controllers in \secref{sec:dp}. A summary of the algorithms is given in \secref{sec:control:summary}.

\subsection{Hierarchical Receding Horizon Controller} \label{sec:hrhc}

In our hierarchical control approach, a high-level path planner finds a progress-optimal, feasible trajectory within a finite number of possibilities and for a horizon of $N$ sampling times. This trajectory is then tracked by an MPC controller employing soft-constraints to ensure feasibility of the low-level optimization at all times. This process is repeated in a receding-horizon fashion, hence we call the controller Hierarchical Receding Horizon Controller (HRHC) in what follows. Details on the individual components are given in the following.

\subsubsection{Path Planning.}
The purpose of the path planner is to generate a trajectory with maximal progress which is feasible for the car's dynamics. It is based on gridding the stationary velocities of the nonlinear system described in Section~\ref{sec:stationaryVelocities}. By solving~\eqref{odeModel} for the stationary lateral velocity points given a set of different longitudinal velocities $\bar{v}_x$ and steering angles $\bar{\delta}$, a library of zero acceleration points is generated. This includes stationary velocities from the normal as well as from the drifting region. Table~\ref{statVeloLibary} shows an excerpt of such a library. The library used in our experiments consists of $N_{\bar{v}} = 95$ stationary points, out of which 26 correspond to drifting equilibria. The stationary points are selected by uniformly gridding the longitudinal velocity $v_x$ between 0.5 and 3.5 m/s, with steps of 0.25 m/s. For each $\bar{v}_x$, about five to nine points are selected in the normal driving region, and up to four drifting points are selected for forward velocities between 1.5 and 2.25 m/s. Note that the library does not change during run-time, hence it can be generated offline.

\begin{table}
\caption{Part of the Path Planning Library}
\centering
\ra{1.3}
\begin{tabular}{l | l l l l}
\toprule
Description & $v_x$ [m/s] & $v_y$ [m/s]& $\omega$ [rad/s]&$\delta$ [rad]\\
\midrule
Left &1.75 &0.1185 & -2.7895 & -0.1000\\
Straight & 1.75 & 0 & 0 & 0\\
Slight Right& 1.75 & -0.0486 & 1.3849 & 0.05\\
Hard Right& 1.75 & -0.3163 & 4.4350 & 0.16\\
Hard Left Drift & 1.75 &  1.1860 & -4.8709 & 0.2\\
\bottomrule
\end{tabular}
\label{statVeloLibary}
\end{table}

\begin{figure}
\centering
\includegraphics[width = 0.8\textwidth]{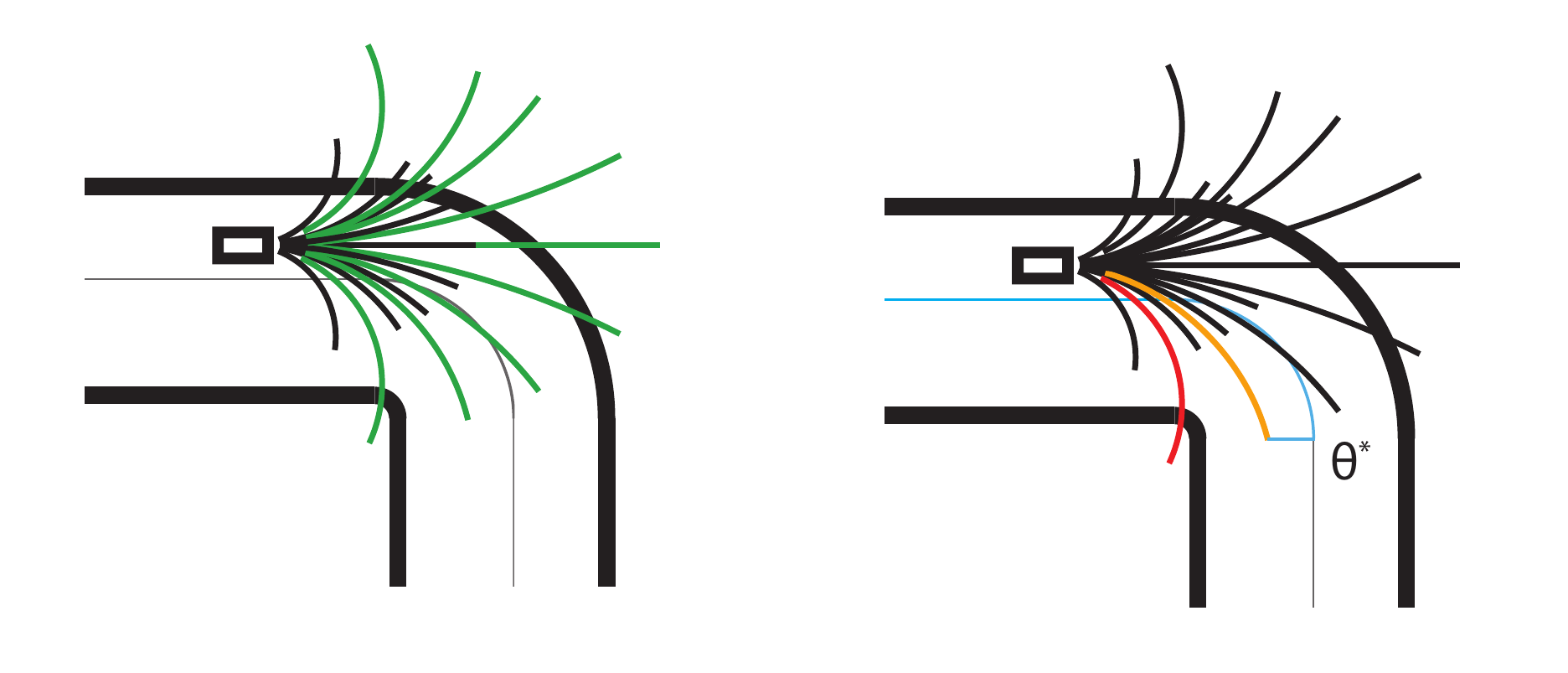}
  \vspace{-1em}
  \caption{Left: Generation of reference trajectories for two different longitudinal velocities (black: slow, green: fast). Right: Solution of the path planner. The red trajectory has the largest progress but leaves the track, so that the orange one is the optimal trajectory maximizing $\theta_N$.}
  \label{pathPlanner1}
\end{figure}

To generate trajectories for the whole horizon, the stationary points are integrated. This is computationally cheap under the assumption that the stationary velocity can be reached within one time step, since in this case all accelerations are zero. This procedure allows us to generate reference trajectories also for unstable drifting points. The integration of the reasonable stationary velocity points leads to a countable set of possible trajectories with a constant turning radius over the horizon. A visualization of such a set of candidate trajectories is shown in Figure~\ref{pathPlanner1}.

We now proceed with finding the best trajectory among the candidates obtained by gridding and integration, which is equivalent to solving the following optimization problem: 
\begin{subequations}\label{ppOptProblem}
\begin{align}
\theta^\star \mydef \max_{j \in 1,\dots,N_{\bar{v}}}\ &  \mathcal{P}\left(X_N^j\,,Y_N^j\right)\,,\label{eq:ppOptProblem:obj} \\
\text{s.t.}\  & x_{0}^j = x\,,\label{eq:ppOptProblem:initial}\\
& x_{k+1}^j = f_{km}( x_k^j, \bar{v}^j) \label{eq:ppOptProblem:model} \,,\,\quad\qquad\qquad k=1,\dots,N \\
&x_{k}^j \in \mathcal{X}_{track}\,,\label{eq:ppOptProblem:trackConstr} \qquad\qquad\qquad\qquad \:\,  k=1,\dots,N\\
&\bar{v}^{j} \in \mathcal{V}(v)\,, \label{eq:ppOptProblem:velConstr}
\end{align}
\end{subequations}
where $x_k \mydef (X_k, Y_k, \varphi_k)$ is the position and orientation at time step $k$ generated by integrating the kinematic part of the bicycle model using the stationary velocity $\bar{v}^j \mydef (\bar{v}_x^j, \bar{v}_y^j, \bar{\omega}^j)$ from the library of stationary points. The integration is denoted by the discrete time version of the model, $f_{km}$ in~\eqref{eq:ppOptProblem:model}. The objective~\eqref{eq:ppOptProblem:obj} is to find the trajectory with the largest progress on the track, which is measured by the projection operator $\mathcal{P}:\reals{2} \rightarrow [0,\,L]$ that calculates the scalar projection of the position $X_k^j\,,Y_k^j$ onto the piecewise linear center line parameterized by the arc length $\theta \in [0,L]$, where $L$ is the length of the center line. All positions of the trajectory $x^j_k,\,k=0,\dots,N,$ need to lie in the set $\mathcal{X}_{track} \subset \reals{2}$ ~\eqref{eq:ppOptProblem:trackConstr} which is the set of all admissible positions inside the track boundaries. Finally, we demand that the stationary velocity $\bar{v}^j$ lies in the set $\mathcal{V} \subset \reals{3}$~\eqref{eq:ppOptProblem:velConstr}, which depends on the current velocity measurement $v\mydef(v_x,v_y,\omega)$ and enforces that the planned velocity is reasonably close to the measured velocity, which helps to deal with the assumption that the stationary velocity can be reached within one time step:
\begin{equation}\label{eq:HRHCtuningConstraints}
\mathcal{V}(v) \mydef \begin{cases} \left\{ (\bar{v}_x, \bar{v}_{y}, \bar{\omega})\in\reals{3}~|~\rho \geq |v_x - \bar{v}_x| \right\} & \text{if}~|\bar{v}_y| \leq \nu \\
 \left\{ (\bar{v}_x, \bar{v}_{y}, \bar{\omega})\in\reals{3}~|~ \sigma_x�\geq |v_x - \bar{v}_x | \,,~\sigma_y \geq |v_y - \bar{v}_y|\,,~\sigma_\omega \geq |\omega - \bar{\omega}| \right\} & \text{otherwise}
\end{cases}\,,
\end{equation}
where the two cases correspond to selecting non-drift or drift trajectories, respectively, based on the lateral steady state velocity $\bar{v}_y$. The tuning parameters $\nu,\rho,\sigma_x,\sigma_y,\sigma_\omega$ are used to determine the behavior of the path planner in terms of selecting candidate trajectories. Choosing parameters that lead to too loose constraints in~\eqref{eq:HRHCtuningConstraints} result typically in physically impossible trajectories, which can result in crashes. Conversely, constraints that are too tight limit the performance of the controller, as the planned trajectories are too conservative.

Problem~\eqref{ppOptProblem} is an integer program with the decision variable $j$ and solved by enumeration in practice by checking the discrete positions against local convex inner approximations of the track set, which is reasonable due to the relatively small number of trajectories (about one hundred in our case). The result of the optimization problem is visualized in Figure~\ref{pathPlanner1}. The path planning allows to efficiently handle drift, as a drifting trajectory is selected if it yields the largest progress.

The path planner has some limitations. Firstly, on one hand, the horizon length has to be quite short, otherwise the planned trajectories become circles and always yield infeasibilities with respect to the track constraints, even if the selected velocity would be well suited for the current state. A horizon which is too short, on the other hand, leads to poor performance as the path planner recognizes a forthcoming curve too late. The second problem is that the whole horizon has the same velocity, which can lead to problems in complicated curve combinations such as chicanes.

\subsubsection{Model Predictive Reference Tracking.} \label{sec:mptc}

\begin{figure}
\centering
\includegraphics[width = 0.8\textwidth]{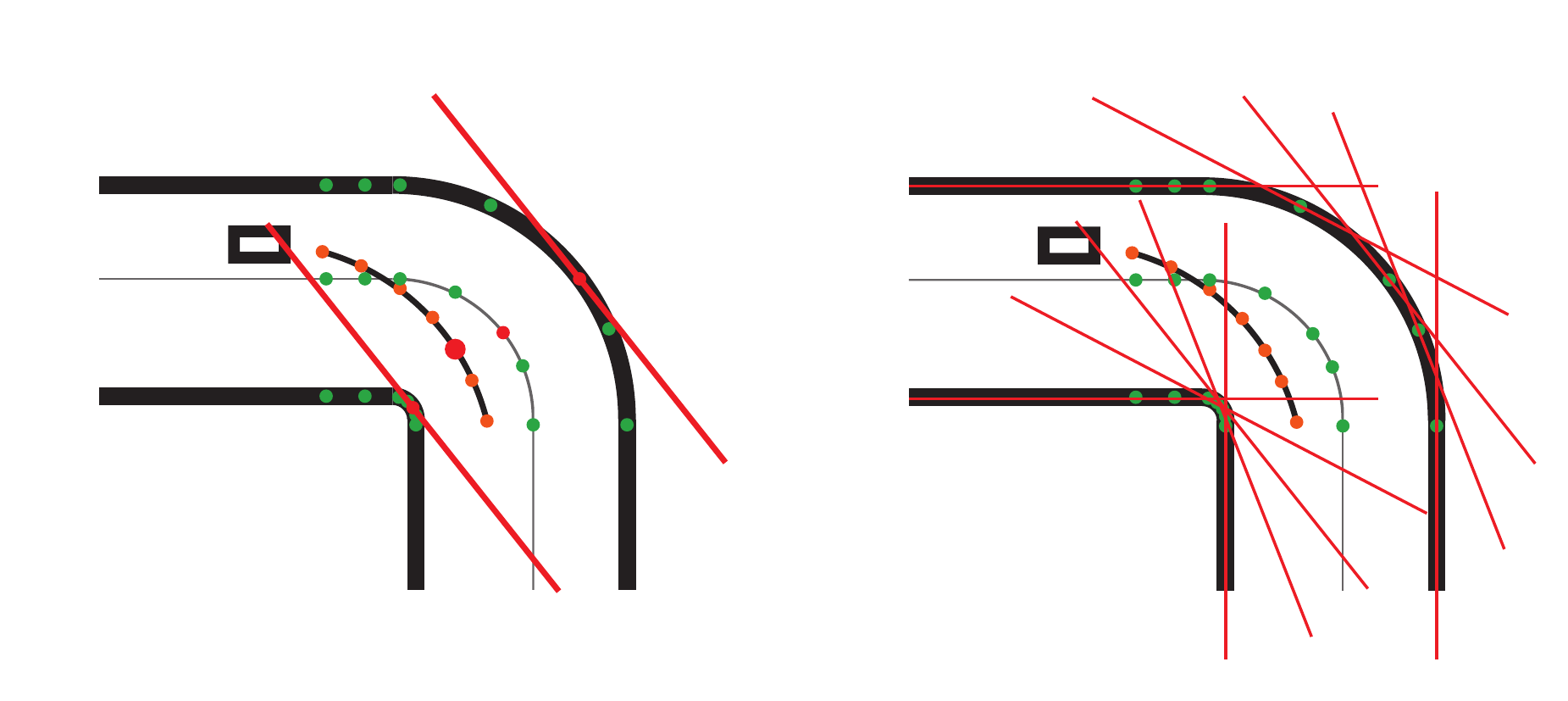}
  \vspace{-1em}
  \caption{Left: Halfspace constraints (red lines) of one point in the horizon~(large red dot). Right: Resulting halfspaces (red lines) of the whole horizon, approximating the track. The green points on the center line correspond to the closest point on the discrete center line, for each point in the horizon. The points on the borders are the projection of the center line points on the borders. }
  \label{halfSpaces}
\end{figure}

The reference trajectory from the path planner is generated under simplifying assumptions, thus it is not possible to directly apply the controls which correspond to the stationary velocity point. In order to efficiently track the reference the following MPC formulation is used, penalizing the deviation from the reference trajectory using a quadratic function:
\begin{subequations}
\label{mptc}
\begin{align}
\min_{x,u,s}\ & \|x_N-x_{N}^{\myref}\|_P^2 + p\|s_N\|_{\infty} + \sum_{k=0}^{N-1} \| x_k-x_{k}^{\myref} \|_Q^2 + q\|s_k\|_{\infty} 
+ \| u_k-u_{k}^{\myref} \|_R^2 \\
\text{s.t.}\ &  x_0 = x\,,\\
&x_{k+1} = A_k x_k + B_k u_k + g_k\,, \qquad \quad \;\;\,  k = 0,...,N-1\label{eq:mptc:dynamics}\\
&F_{k} x_k \leq f_{k} + s_k\,, \qquad \qquad  \qquad  \qquad \;  k = 1,...,N\label{eq:mptc:softBorders}\\
&s_k \geq 0\,,\qquad \qquad \qquad  \qquad  \qquad \qquad   k = 1,...,N\\
&\underline{x} \leq x_k \leq \bar{x} \,, \qquad \qquad \qquad  \qquad \qquad   k = 1,...,N \label{eq:mptc:otherconstraints}\\
&\underline{u} \leq u_k \leq \bar{u}\,, \qquad \qquad \qquad \qquad \qquad k = 0,...,N-1\label{eq:mptc:otherconstraints2}
\end{align}
\end{subequations}
where $Q\in\pdmat{6}$, $R\in\pdmat{2}$, $P\in\pdmat{6}$ are tuning matrices and $p,q\in\posreals{}$ are penalties for the soft constraints~\eqref{eq:mptc:softBorders}, which are discussed below. To capture as much of the nonlinear model as possible while maintaining computational tractability, the model is linearized around the reference trajectory, cf.~\eqref{eq:mptc:dynamics}. It is necessary to linearize around a trajectory rather than a single operating point, since both position and orientation can vary significantly over the horizon. The reference trajectory is by construction feasible with respect to the track. In order to also guarantee that the trajectory of the optimal MPC solution satisfies track constraints, the $X$ and $Y$ states are limited to stay inside the track if possible. This is achieved by limiting each point in the horizon to lie within two parallel half spaces~\eqref{eq:mptc:softBorders}, which leads to two affine inequality constraints per time step (one for the right and one for the left border), resulting in a convex feasible set. Figure~\ref{halfSpaces} depicts these border constraints~\eqref{eq:mptc:softBorders}, which we formulate as soft constraints in order to avoid infeasibility problems in practice. By using an exact penalty function (an infinity norm in our case), we ensure that the resulting optimization problem is always feasible, and that the original solution of the hard constrained problem is recovered in case it would admit a solution~\cite{Luenberger2003,Kerrigan2000}. Finally, the control inputs are limited within their physical constraints, and all other states are limited within reasonable bounds to avoid convergence problems of the solver due to free variables~\eqref{eq:mptc:otherconstraints},~\eqref{eq:mptc:otherconstraints2}. 

Problem~\eqref{mptc} can be reformulated as a convex QP and thus efficiently solved. In this work FORCES~\cite{Domahidi2012,FORCEScodegen}, is employed to generate tailored C-code that solves instances of~\eqref{mptc} where parameters are the initial state $x$, the matrices defining the dynamics along the reference trajectory ($A_k,B_k,g_k$) as well as the changing halfspace constraints $F_k,f_k$.

\subsection{Model Predictive Contouring Control} \label{sec:mpcc}

Contouring control is used in various industrial applications such as machine tools for milling and turning~\cite{Koren1997} or laser profiling~\cite{Ko1999}. In these applications, the challenge is to compute inputs that control the movement of the tool along a reference path. The latter is given only in spatial coordinates; the associated velocities, angles, etc. are calculated and imposed by the control algorithm. This is different from tracking controllers in that the controller has more freedom to determine the state trajectories to follow the given path, for example to schedule the velocity, which in tracking is defined by the reference trajectory.

The contour following problem can be formulated in a predictive control framework to incorporate constraints, see e.g.~\cite{Faulwasser2009}. In this work, we adapt the particular formulation of the model predictive contouring control (MPCC) framework from~\cite{Lam2010} to obtain a high-performance controller for autonomous racing, which maximizes the travelled distance on the reference path within the prediction horizon. In our experiments, we use the center line as a reference path, but employ it merely as a measure of progress by selecting low weights on the tracking (contouring) error. As a result, the driven trajectory is very similar to those driven by expert drivers. The advantage of this approach is that path planning and path tracking can be combined into one nonlinear optimization problem, which can be solved in real-time by approximating the NLP using local convex QP approximations at each sampling time~\cite{Diehl2002b}. The resulting QPs can be formulated with multistage structure, which is exploited by FORCES~\cite{Domahidi2012,FORCEScodegen} to obtain solution times in the range of a few milliseconds.

Before posing the MPCC problem, a few preliminaries are introduced such as the parameterization of the reference path and the definition of useful error measures. 

\subsubsection{Parameterization of Reference Trajectory.}

The reference path is parameterized by its arc length $\theta \in [0,\,L]$ using third order spline polynomials, where $L$ is the total length. The splines are obtained by an offline fitting of the center line. Using this parameterization, we can obtain any point $X^\myref(\theta),\,Y^\myref(\theta)$ on the center line by evaluating a third order polynomial for its argument~$\theta$. The angle of the tangent to the path at the reference point with respect to the $X$-axis,
\begin{equation}
\Phi(\theta)\mydef \arctan \left\{ \frac{\partial Y^\myref(\theta)}{\partial X^\myref(\theta)} \right\}\,,
\label{eq:PhiOfTheta}
\end{equation}
is also readily available. This parameterization leads to an accurate interpolation within the known points of the reference path, and it is more accurate than the piece-wise linear parameterization employed for the HRHC in \secref{sec:hrhc}, which only uses a linear interpolation between the points.

\subsubsection{Error measures.} \label{sec:mpcc:errorMeasures}

\begin{figure}
\centering
\includegraphics[width = 0.8\textwidth]{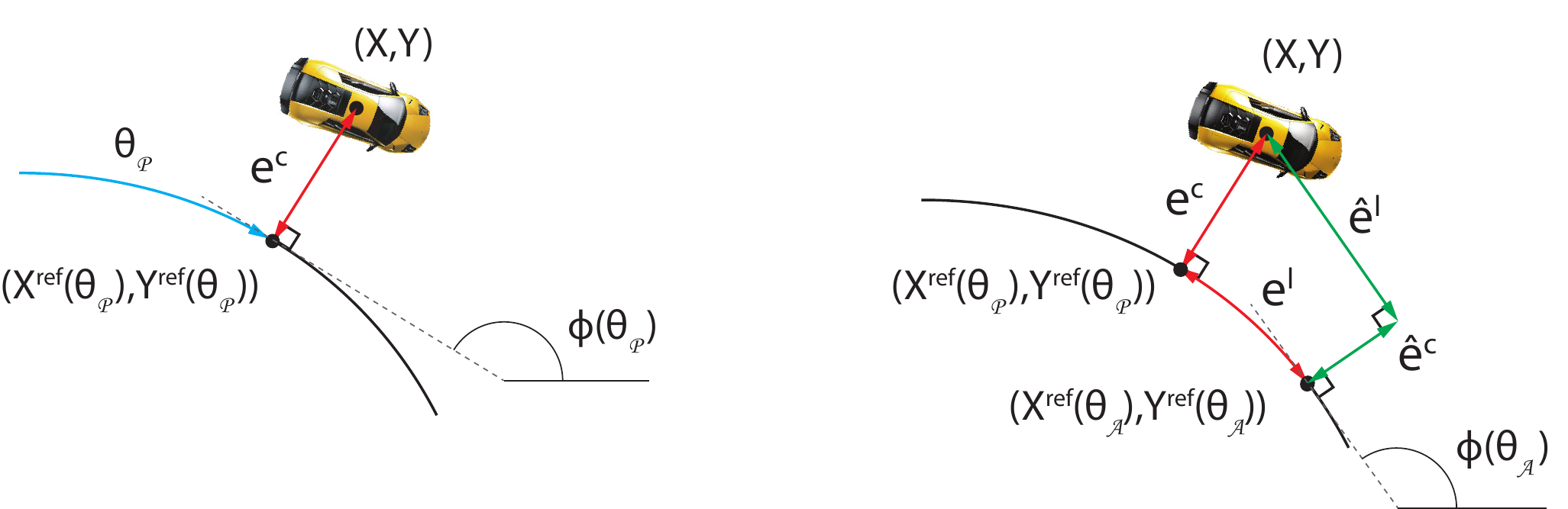}
  \caption{Contouring error $e^c$ (left) and lag error $e^l$ (right) with linear approximations $\hat{e}^c$ and $\hat{e}^l$. }
  \label{mpccError1}
\end{figure}

In order to formulate the MPCC problem, error measures are needed that define the deviation of the car's current position $X,\,Y$ from the desired reference point $X^\myref(\theta),\,Y^\myref(\theta)$. We use the same definitions as in~\cite{Lam2010}, but give another derivation here. Let $\mathcal{P}:\reals{2} \rightarrow [0,\,L]$ be a projection operator on the reference trajectory defined by
\begin{equation}
\mathcal{P}(X,Y) \mydef \arg \min_\theta (X-X^\myref(\theta))^2 + (Y-Y^\myref(\theta))^2\,.
\label{eq:proj}
\end{equation}
For brevity, we define $\theta_{\mathcal{P}} \mydef \mathcal{P}(X,Y)$. The orthogonal distance of the car from the reference path is then given by the \emph{contouring error}
\begin{equation}
e^c(X,Y,\theta_{\mathcal{P}}) \mydef \sin\left(\Phi(\theta_{\mathcal{P}})\right)\left(X - X^\myref(\theta_{\mathcal{P}})\right) 
- \cos\left(\Phi(\theta_{\mathcal{P}})\right)\left(Y - Y^\myref(\theta_{\mathcal{P}})\right)\,,
\label{eq:contErr}
\end{equation}
where $\Phi(\cdot)$ is defined in~\eqref{eq:PhiOfTheta}. The contouring error is depicted in the left picture in Figure~\ref{mpccError1}. 

The projection operator~\eqref{eq:proj} is not well suited for use within online optimization algorithms, as it resembles an optimization problem itself. Thus an approximation $\theta_{\mathcal{A}}$ of $\theta_{\mathcal{P}}$ is introduced, which is an independent variable determined by the controller. 
For this approximation to be useful, it is necessary to link $\theta_{\mathcal{A}}$ to $\theta_{\mathcal{P}}$ via the \emph{lag error}
\begin{equation}
e^l(X,Y,\theta_{\mathcal{A}}) \mydef | \theta_{\mathcal{A}} - \theta_{\mathcal{P}}|\,, \label{eq:lagErr}
\end{equation}
which measures the quality of the approximation. The lag error is depicted in the right picture in Figure~\ref{mpccError1} (red segment on reference path). 

In order to be independent of the projection operator~\eqref{eq:proj}, the contouring \eqref{eq:contErr} and the lag error \eqref{eq:lagErr} can be approximated as a function of the position $X,\,Y$ and the approximate projection $\theta_{\mathcal{A}}$, which are variables that can be controlled by the MPCC controller:
\begin{subequations}\label{eq:mpcc:approximateErrorMeasures}
\begin{align}
e^c &\approx \hat{e}^c(X,Y,\theta_{\mathcal{A}}) \mydef \sin\left(\Phi(\theta_{\mathcal{A}})\right)\left(X - X^\myref(\theta_{\mathcal{A}})\right) - \cos\left(\Phi(\theta_{\mathcal{A}})\right)\left(Y - Y^\myref(\theta_{\mathcal{A}})\right)\,,\label{eq:appContErr}\\
e^l &\approx \hat{e}^l(X,Y,\theta_{\mathcal{A}}) \mydef -\cos\left(\Phi(\theta_{\mathcal{A}})\right)\left(X - X^\myref(\theta_{\mathcal{A}})\right) - \sin\left(\Phi(\theta_{\mathcal{A}})\right)\left(Y - Y^\myref(\theta_{\mathcal{A}})\right)\,.\label{eq:appLagErr}
\end{align}
\end{subequations}
The approximate contouring error $\hat{e}^c$ \eqref{eq:appContErr} and the approximate lag error $\hat{e}^l$ \eqref{eq:appLagErr} are defined as the orthogonal and tangential component of the error between $X^\myref(\theta_{\mathcal{A}}),\,Y^\myref(\theta_{\mathcal{A}})$ and the position  $X,\,Y,$ see the right picture in Figure \ref{mpccError1}.

\subsubsection{MPCC problem.} With these error measures in place, we now formulate the model predictive contouring control problem for autonomous racing. The objective is a trade-off between the quality of path-following (low contouring error) and the amount of progress achieved over a finite horizon of $N$ sampling times, subject to model dynamics, track- and input constraints:
\begin{subequations}\label{eq:mpcc}
\begin{align}
\min\ & \sum_{k=1}^{N} \left\{ \|e_k^c(X_k,Y_k,\theta_{\mathcal{P}})\|^2_{q_c}\right\} - \gamma \theta_{\mathcal{P},N} \label{eq:mpcc:objective}\\
\text{s.t.}\ & x_0 = x\,,\label{eq:mpcc:initialState} \\ 
& x_{k+1} = f(x_k,u_k)\,, \qquad \quad k = 0,...,N-1 \label{eq:mpcc:model}\\
& F_k x_k \leq f_k \,, \qquad \qquad \quad \;\;\; k = 1,...,N\label{eq:mpcc:borderConstr} \\
&\underline{x} \leq x_k \leq \bar{x} \,, \qquad  \qquad \quad \;\; k = 1,...,N\label{eq:mpcc:otherconstraints}\\
&\underline{u} \leq u_k \leq \bar{u}\,, \qquad \qquad \quad \;\; k = 0,...,N-1 \label{eq:mpcc:otherconstraints2}
\end{align}
\end{subequations}
where $X_k\,,Y_k$ is the position of the car at time step $k$ determined by the nonlinear model $f$ in~\eqref{eq:mpcc:model}, which is the discrete-time version of \eqref{odeModel} with piece-wise constant control inputs. The contouring error $e_k^c(X_k,Y_k,\theta_{\mathcal{P}})$ in the objective~\eqref{eq:mpcc:objective} is defined in~\eqref{eq:contErr}, and $\theta_{\mathcal{P},k}$ is the associated path parameter such that $X^\myref(\theta_{\mathcal{P},k}),\,Y^\myref(\theta_{\mathcal{P},k})$ is the orthogonal projection of $X_k,\,Y_k$ onto the reference path. The two objectives (maximum progress and tight path following) are traded off by the weights $\gamma\in\posreals{}$ and $q_c\in\posreals{}$, respectively. Constraints \eqref{eq:mpcc:borderConstr} are the parallel half space constraints for containing the position in the allowed corridor as described in \secref{sec:mptc}, while \eqref{eq:mpcc:otherconstraints}, \eqref{eq:mpcc:otherconstraints2} corresponds to~\eqref{eq:mptc:otherconstraints},  \eqref{eq:mptc:otherconstraints2}, limiting states and inputs to physically admissible values. 

Due to the implicit dependency of the objective~\eqref{eq:mpcc:objective} on the projection operator~\eqref{eq:proj}, optimization problem~\eqref{eq:mpcc} is a bi-level NLP, which is too complex to solve in real-time. Thus the idea is to use the approximate projection $\theta_{\mathcal{A},k}$ instead of $\theta_{\mathcal{P},k}$ as introduced in \secref{sec:mpcc:errorMeasures}, and to control the approximation quality by adding a cost on the lag error~\eqref{eq:lagErr} to the objective. In order to allow for forming the lag error at each time step in the prediction horizon, it is necessary to introduce an integrator state with dynamics $\theta_{\mathcal{A},k+1} = \theta_{\mathcal{A},k} + v_k/T_s$, where $v_k$ can be interpreted as the \emph{projected velocity}, and $\theta_{\mathcal{A},k}$ as the \emph{state of progress} at time $k$, respectively. This approximation reduces~\eqref{eq:mpcc} to an optimal control problem in form of an NLP that is amenable for a real-time implementation: 
\begin{subequations} \label{eq:nlocp}
\begin{align}
\min \ & \sum_{k=1}^{N} \| \hat{e}^c_k(X_k,Y_k,\theta_{\mathcal{A},k}) \|_{q_c}^2 + \| \hat{e}^l_k(X_k,Y_k,\theta_{\mathcal{A},k}) \|_{q_l}^2 - \gamma v_k T_s + \| \Delta u_k \|_{R_u}^2 + \| \Delta v_k \|_{R_v}^2 -\gamma v_0 T_s \label{eq:mpcc:approx:obj}\\
\text{s.t.}\ & \eqref{eq:mpcc:initialState}\,, \quad \theta_0 = \theta\,, \label{eq:mpcc:approx:initialStates}\\
&\eqref{eq:mpcc:model}\,, \quad \theta_{\mathcal{A},k+1} = \theta_{\mathcal{A},k} + \frac{v_k}{T_s}\,, \qquad \qquad k=0,...,N-1\label{eq:mpcc:approx:dynamics}\\
& \eqref{eq:mpcc:borderConstr}\,,\quad\eqref{eq:mpcc:otherconstraints}\,, \quad 0\leq \theta_k \leq L\,, \qquad\qquad  k = 1,...,N\label{eq:mpcc:approx:constraints1}\\
&\eqref{eq:mpcc:otherconstraints2}\,, \quad 0 \leq v_k \leq \bar{v}\,,\qquad\qquad \qquad\qquad k = 0,...,N\label{eq:mpcc:approx:constraints2}
\end{align}
\end{subequations}
where $\Delta u_k \mydef u_k - u_{k-1}$ and $\Delta v_k \mydef v_k - v_{k-1}$. Note that the approximate contouring error $\hat{e}^c_k(X_k,Y_k,\theta_{\mathcal{A},k})$ \eqref{eq:appContErr} is used in the objective~\eqref{eq:mpcc:approx:obj}. Furthermore, we have replaced the maximization of the final progress measure, $\theta_{\mathcal{P},N}$, by $\sum_{k=0}^{N-1} v_k T_s$, which is equivalent if the approximation is accurate. Sensible lower and upper bounds on $\theta_k$ and $v_k$ are imposed to avoid spurious solutions of the NLP, with $\bar{v}$ denoting the largest possible progress per sampling time. The cost on the lag error $\hat{e}^l_k(X_k,Y_k,\theta_{\mathcal{A},k})$ in~\eqref{eq:mpcc:approx:obj} links the state of progress to the dynamics of the car. To ensure an accurate progress approximation and thus a strong coupling between the cost function and the car model, the weight on the lag error $q_l \in \posreals{}$ is chosen high as suggested in~\cite{Lam2010}. Furthermore, a cost term on the rate of change of the inputs is added to the objective~\eqref{eq:mpcc:approx:obj} in order to penalize fast changing controls, which helps to obtain smooth control inputs in order to prevent amplifying unmodeled dynamics.

\subsubsection{Solving the MPCC problem.} \label{sec:mpcc:solving} 
In order to solve the nonlinear optimal control problem~\eqref{eq:nlocp} in real-time, local convex approximations of~\eqref{eq:nlocp} in form of the following QPs are built at each sampling time by linearization of nonlinear terms:
\begin{subequations} \label{qpForm}
\begin{align}
\min_{x,u,\theta,v,s}\ & \sum_{k=1}^{N} \begin{bmatrix} x_k \\ \theta_{\mathcal{A},k} \end{bmatrix}^T \Gamma_k \begin{bmatrix} x_k \\
\theta_{\mathcal{A},k} \end{bmatrix} + c_k^T \begin{bmatrix} x_k \\ \theta_{\mathcal{A},k} \end{bmatrix} - \gamma v_{k} T_s +  \begin{bmatrix} \Delta u_k \\
\Delta v_k \end{bmatrix}^T R \begin{bmatrix} \Delta u_k \\ \Delta v_k \end{bmatrix} + q\|s_k\|_{\infty} - \gamma v_0 T_s \label{eq:mpcc:qp:obj}\\
\text{s.t.}\ & x_0 = x\,, \quad \theta_{\mathcal{A},0} = \theta\,,\\
&x_{k+1} = A_k x_k + B_k u_k + g_k\,, \qquad \qquad\qquad\; k=0,...,N-1\\
& \theta_{\mathcal{A},k+1} = \theta_{\mathcal{A},k} + \frac{v_k}{T_s}\,,\qquad\qquad\qquad\qquad\quad\; k=0,...,N-1\\
&F_{k} x_k \leq f_{k} + s_k\,,\qquad\qquad\qquad\qquad\qquad\quad k=1,...,N\label{eq:mpcc:borderConstrSoft}\\
&s_k \geq 0\,, \qquad\qquad\qquad\qquad\qquad\qquad\qquad\;\;\, k=1,...,N\\
&\underline{x} \leq x_k \leq \bar{x} \,, \quad 0\leq \theta_{\mathcal{A},k} \leq L \,, \qquad\qquad\quad\;\; k=1,...,N\\
&\underline{u} \leq u_k \leq \bar{u}\,, \quad 0 \leq v_k \leq \bar{v} \,, \qquad\qquad\qquad\;\; k=0,...,N
\end{align}
\end{subequations}
where $\Gamma_k \in \psdmat{7}$ is formed by the quadratic part of the linearized contouring and lag error cost function from~\eqref{eq:mpcc:approx:obj} and $c_k \in \reals{7}$ stems from the linear part, respectively. In order to keep the linearization error small, we use an LTV approximation of the dynamics~\eqref{odeModel} as well as of the contouring and lag errors~\eqref{eq:mpcc:approximateErrorMeasures}. Each nonlinear function is linearized around the output of the last QP iteration shifted by one stage. The measurement $x=x_0$ is used as the first linearization point, and the last input of the previous iteration is kept constant to generate a new last input. The linearization point for the terminal state $x_N$ is calculated by simulating the nonlinear model for one time step. The track or corridor constraints~\eqref{eq:mpcc:borderConstrSoft} are formulated by two half space constraints tangential to the track per time step, as described in Section~\ref{sec:hrhc}. These constraints are formulated as soft constraints, with slack variables $s_k \in \reals{2}$ and a corresponding infinity-norm penalty in the objective~\eqref{eq:mpcc:qp:obj} weighted by $q \in \posreals{}$, which is chosen quite high to recover the behavior of the hard constrained problem whenever possible~\cite{Luenberger2003,Kerrigan2000}. 

Note that the described re-linearization scheme is essentially the basic real-time iteration from~\cite{Diehl2002b}, where a nonlinear continuous-time optimal control problem is solved using a multiple shooting technique, with an exponential map as integrator, and in every time step only one step of a sequential QP (SQP) scheme is performed. Hence the convergence properties that are known for the real-time iteration~\cite{Diehl2005} hold also for the MPCC problem~\eqref{eq:nlocp}.

The local QP~\eqref{qpForm} is solved using FORCES~\cite{Domahidi2012,FORCEScodegen}. Due to the rate costs, the problem has to be lifted to obtain the multistage structure, for which FORCES has been designed, by introducing a copy of the previous control input at each stage.

\subsection{Obstacle Avoidance} \label{sec:dp}

\begin{figure}
\centering
\includegraphics[width = \textwidth]{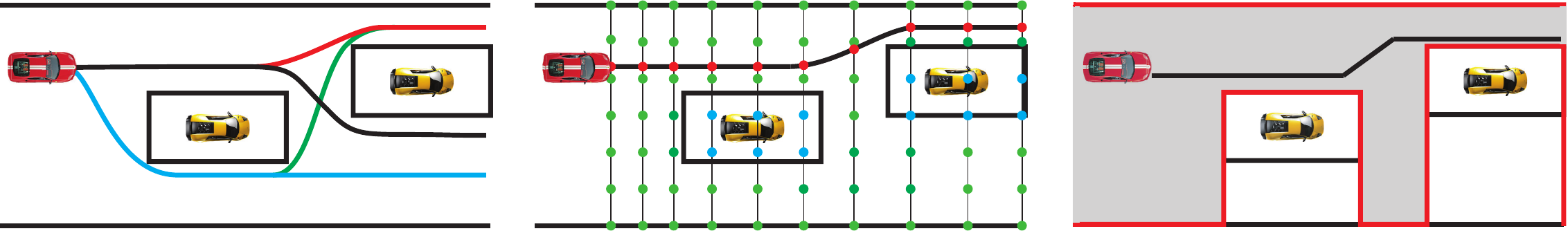}
  \caption{Left: Combinatorial nature of the overtaking problem, as it is possible to overtake each opponent on the left or the right side. Center: Spatial-temporal grid for the shortest path problem solved by dynamic programming, where it is not allowed to visit the blue grid points. Right: Optimal path of the DP (in black) and the resulting corridor, which is issued to the HRHC (\secref{sec:hrhc}) or MPCC (\secref{sec:mpcc}).}
  \label{DP}
\end{figure}

The two controllers presented in this section plan a path within the track, or a corridor, by taking into account constraints \eqref{eq:ppOptProblem:trackConstr}, \eqref{eq:mptc:otherconstraints} or \eqref{eq:mpcc:borderConstrSoft} on the position of the car.  Thus it is possible to include obstacle avoidance by adapting this corridor online depending on the current constellation of opponents. Generating the optimal corridor based on the position of the obstacles is in general hard, because of the fundamental problem that the decision on which side to overtake is non-convex and of combinatorial nature if more than one opponent is involved. This situation is depicted in Figure~\ref{DP} in the left picture.

In this work, we employ a high level path planner based on dynamic programming (DP) to decide on which side to overtake the obstacles. Based on the result of the DP, the adapted corridor is given to the controllers of \secref{sec:hrhc} or \secref{sec:mpcc}. The high level path planner solves a shortest path problem~\cite{Bertsekas1995} on a spatial-temporal grid, cf. Figure \ref{DP}. To incorporate the planned path of the lower level controller into the DP, the grid is built up based on the last output of the path planner in the case of the HRHC or the last QP solution in the case of the MPCC. The DP minimizes the travelled distance, and additionally the deviation from the last plan of the lower level controller, such that the DP path is close to the state prediction of the lower level controller. This ensures that the linearizations stay approximately valid most of the time. Based on the optimal trajectory from the DP, the corridor constraints can easily be identified, see the right picture in Figure \ref{DP}.

Due to space restrictions, we do not go further into detail, but many references exist that solve the corridor problem with different approaches. For example, a similar problem had to be solved during the DARPA Urban challenge where the goal was autonomous driving in city traffic. In \cite{Buehler2009} several approaches to this problem are presented, for example based on model predictive trajectory generation algorithms~\cite{Ferguson_2008}, modified $A^\star$-algorithms~\cite{Montemerlo_2008}~or RRTs~\cite{Kuwata_2008}.

\subsection{Summary} \label{sec:control:summary}

In this main part of the paper, we have presented two approaches for automatic racing. Both approaches are model-based and maximize progress within a given time horizon. They directly incorporate track and obstacle constraints by using two parallel affine inequalities (slabs) representing the feasible set of positions at each time. The major difference between the approaches is that the first (cf.~\secref{sec:hrhc} and~Algorithm~\ref{alg:HRHC}) is a two-level approach with a separate path planning mechanism combined with a reference tracking MPC controller, while the second approach (cf. \secref{sec:mpcc} and Algorithm~\ref{alg:MPCC}) formulates both tasks (path planning and path tracking) into one nonlinear optimal control problem based on model predictive contouring control.

\begin{algorithm}
\begin{subalgorithm}{0.48\textwidth}	
		\captionof{algorithm}{HRHC}
		\label{alg:HRHC}
		\begin{algorithmic}[1]
   				\State get current position and velocities \label{alog:hrhc:beginning}
   				\State compute delay compensation \label{alog:hrhc:comp}
   				\Function{BorderAdjustment}{}
   					\State shift old planned path by one stage
   					\State get position of opposing cars
   					\State run DP and get new borders (Sec.~\ref{sec:dp})
   				\EndFunction
   				\State run path planner~\eqref{ppOptProblem}
   				\State linearize and discretize the model, i.e. build problem~\eqref{mptc}
   				\State solve QP using FORCES~\eqref{mptc}
   				\State send first control at the end of time slot
   				\State \textbf{go to} step \ref{alog:hrhc:beginning}
		\end{algorithmic}
\end{subalgorithm}
\hfill
\begin{subalgorithm}{0.48\textwidth}
		\captionof{algorithm}{MPCC}
		\label{alg:MPCC}
		\begin{algorithmic}[1]
   				\State get current position and velocities \label{alog:mpcc:beginning}
   				\State compute delay compensation  \label{alog:mpcc:comp}
   				\State augment old QP output (Section~\ref{sec:mpcc:solving})
   				\Function{BorderAdjustment}{}
   					\State get position of opposing cars
   					\State run DP and get new borders (Sec.~\ref{sec:dp})
   				\EndFunction
   				\State linearize and discretize the model and the cost function, i.e. build problem~\eqref{qpForm}
   				\State solve QP using FORCES~\eqref{qpForm}
   				\State send first control at the end of time slot
   				\State \textbf{go to} step \ref{alog:mpcc:beginning}
		\end{algorithmic}
\end{subalgorithm}
\captionsetup{labelformat=alglabel}
\caption{HRHC and MPCC algorithm}
\end{algorithm}

Both controllers send the new control input at the end of the sampling period. In order to compensate for this fixed delay, the nonlinear model is simulated forward for this time interval at the beginning of the algorithm, see line \ref{alog:hrhc:comp} of Algorithm \ref{alg:HRHC} and \ref{alg:MPCC}.  The nonlinear model is integrated using a second-order Runge-Kutta method with the current state estimate from the Kalman filter as the initial condition.

\section{Results} \label{sec:results}
In this section the performance of both controllers is evaluated on an experimental testbed using 1:43 scale RC race cars. Details on the experimental setup, the particular implementation and closed-loop performance are given in the following.

\subsection{Experimental Setup}
For our experiments, we use the Kyosho \dnano~cars, which are quite sophisticated with a front and rear suspension and a rear axle differential. The cars are able to reach speeds of over 3\,m/s on straights, which corresponds to an upscaled speed of about 465\,km/h. The wide speed range, the high maximal forward velocity and the small scale introduce further complexity which is not present in full size testbeds.

\begin{figure}
\centering
\includegraphics[width = 0.8\textwidth]{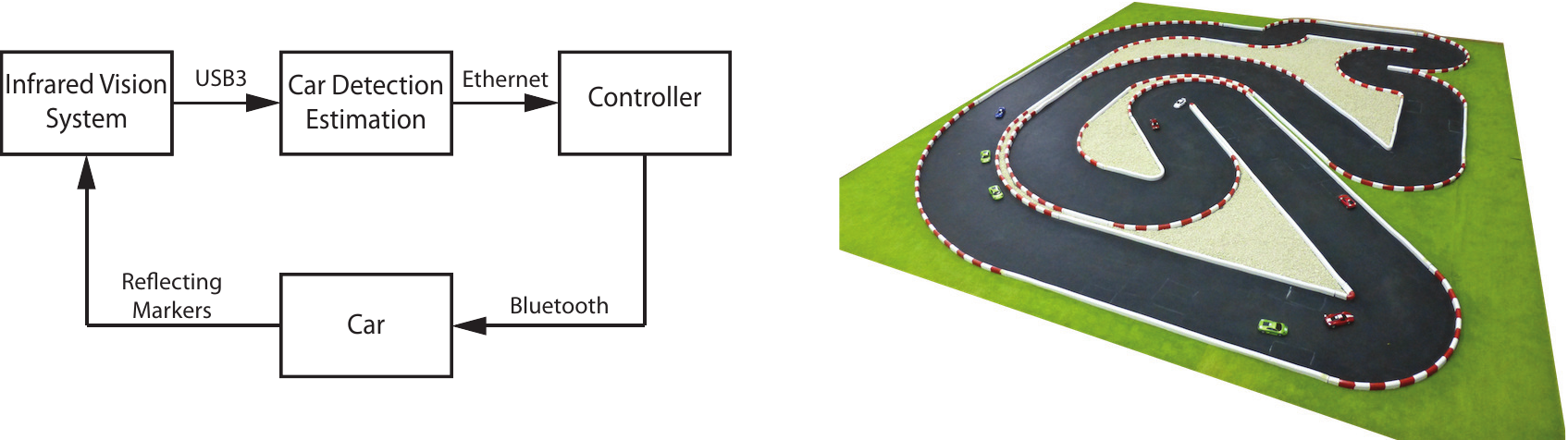}
  \caption{Visualization of the control loop and a picture of the used track}
  \label{ControlLoop}
\end{figure}

The testbed further consists of an infrared camera tracking system, a control board and a custom built race track of length $L =18.43$\,m (measured at the center line), see Figure~\ref{ControlLoop}. A PointGrey Flea3 camera captures 100\,frames per second with an accuracy of below 4\,mm. A wide angle lens with an infrared filter is used to capture the 4 by 4 meter track. Reflecting markers are arranged on the cars in different unique patterns and used to identify the car and to measure the position and angle of the cars. An Extended Kalman Filter (EKF) is used for state estimation, filtering the vision data and estimating the longitudinal and lateral velocity as well as the yaw rate. The vision data of all cars is passed over an Ethernet connection to a control PC, which runs one of the two control algorithms described in \secref{sec:control} together with the corridor planner solving the DP problem. The calculated inputs are then sent via Bluetooth to the car.  As the original electronics of the \dnano~cars does not feature an open digital communication link, we have replaced it by a custom made PCB, featuring a Bluetooth chip, current and voltage sensors, H-bridges for DC motor actuation, and an ARM Cortex M4 microcontroller for the low level control loops such as steering servo and traction motor control.

\subsection{Implementation}

The HRHC is implemented in ANSI C on an embedded platform using an ARM A9 chip (Exynos 4412) at 1.7 GHz. The chip is identical to the one used in Samsung Galaxy S3 smartphones. The runtime environment was Ubuntu Linux version 12.11, and the code was compiled with the GCC compiler 4.6.3 with option \texttt{-O2}. The sampling time of the controller is $T_s=20$\,ms, while the discretization time of the path planner and the tracking MPC controller is 25\,ms. The horizon length is $N=14$ time steps, which gives an ahead prediction of $0.35$\,s. The tracking MPC controller results in a multistage QP with $9$ variables per stage, i.e. $132$ variables in total. We use a library of $95$ candidate trajectories, out of which $26$ are drifting stationary velocity points.

The MPCC does not have the limitation of the HRHC path planner, or in other words, the acceleration is not zero over the whole horizon, thus it is possible to use longer horizons which should increase the performance. Hence the controller runs with a horizon length of $N=40$ with a sampling time of $T_s=20$\,ms, which corresponds to a $0.8$\,s ahead prediction. With the re-linearization method presented in \secref{sec:mpcc:solving}, the discretization and sampling time have to be the same, as the re-linearization is based on the last QP solution. The resulting QP has 14 variables per stage and 570 in total.

The MPCC is implemented on a newer gereration embedded platform using an Exynos 5410 chip based on the ARM Cortex A15 architecture and running at 1.6 GHz. The chip is identical to the one used in the Korean version of the Samsung Galaxy S4 smartphone. The embedded board runs Ubuntu Linux version 12.04. The ANSI C code of the MPCC controller has been compiled with GCC 4.6.3 with option~\texttt{-O3}\,. To improve the convergence of the QP, the problem data is scaled to lie between $-1$ and $+1$, and the objective (Hessian and linear term) is scaled by a factor of $0.1$ to reduce large entries occurring due to the lag error term. Note that these scalings do not change the minimizer. The optimality criteria of the FORCES solver are adjusted such as to terminate once a sufficiently good solution has been found, with the residuals of the equality and inequality constraints set to $10^{-5}$, and with the duality gap set to $10^{-4}$.

\subsection{Single Car Racing}

\begin{figure}
\centering
\centering
\footnotesize
HRHC controller (\secref{sec:hrhc}) \hspace{3.75cm} MPCC controller (\secref{sec:mpcc})\\
\normalsize
\includegraphics[width = \textwidth, trim= 0cm 0cm 0cm 1.2cm, clip]{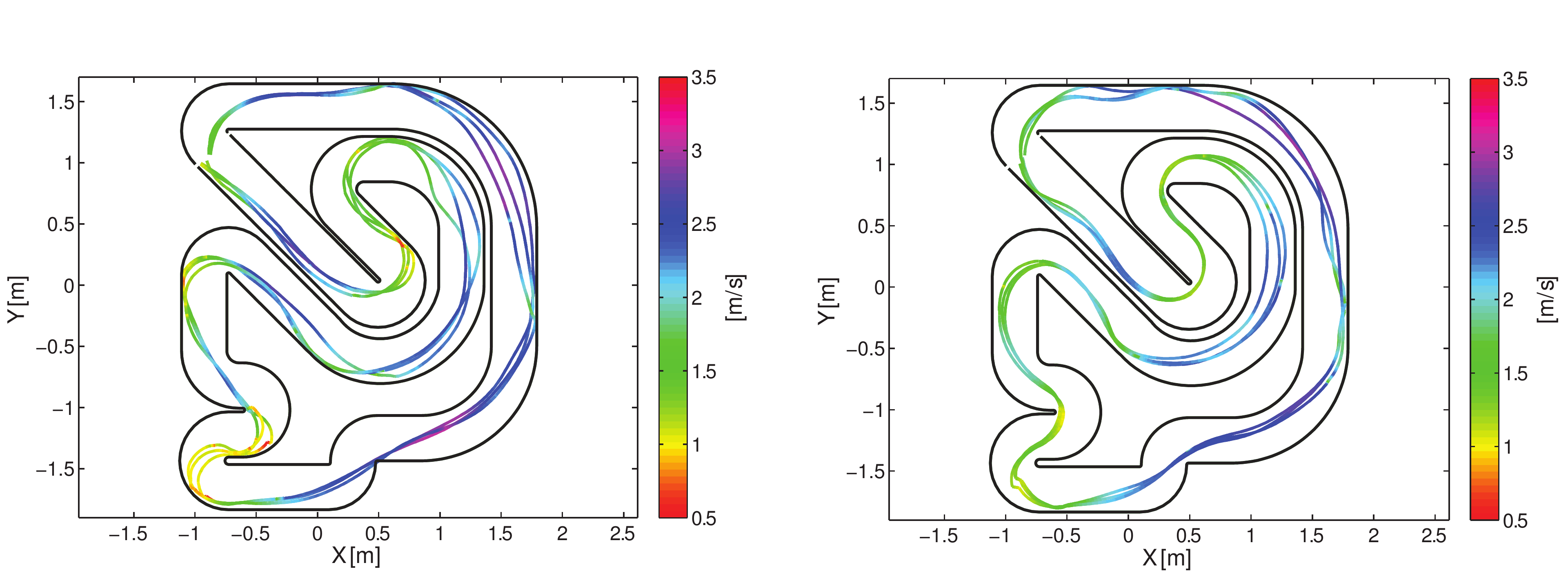}
  \caption{The driven trajectory with velocity profile for 3 different laps. A video of the HRHC controller is available at \url{https://www.youtube.com/watch?v=ioKTyc9bG4c}.}
  \label{obFree}
\end{figure}

The driven trajectories for a single car are depicted in Figure~\ref{obFree}, with the velocity profile encoded in colors. Depicted are three laps, for which the HRHC achieves a lap time between 9.5 to 9.8\,s while the MPCC yields lap times between 8.9 and 9.2\,s.

To better understand the driven trajectory of the HRHC, the path planner of the HRHC has to be analyzed. The path planner has a comparably short prediction horizon, and all velocities over the whole horizon are constant, which leads to turns of constant radius. Such a path planner works fine for 90$^{\circ}$ and even up to 180$^{\circ}$  curves. However, the car tends to drive to the outer border after the curve, as this allows driving a wider radius at a faster longitudinal velocity and thus maximizes the progress. However, this is a limitation if the car should drive a combination of curves, where the position of the car at the end of the first curve is essential for a low overall time. This problem can be seen in Figure~\ref{obFree}, for example in the chicane in the lower left corner. Furthermore, due to the short lookahead, the HRHC path planner is not able to prevent such an ill positioning relative to a curve ahead. Due to the suboptimal position relative to the curve, the controller has to brake more compared to a car on the ideal line and it is even possible that the car touches the borders, see the narrow 180$^{\circ}$  curve in the center of the track in Figure \ref{obAvoid}.

These limitations are not present in the MPCC approach, where a comparably long horizon is employed. Unlike the HRHC, the velocity is not fixed over the horizon but computed for each time step as the result of the NLP. As a consequence of these features, the MPCC is able to plan a path even through a complicated combination of curves. This results in a trajectory which is closer to an ideal line in a least curvature sense, especially through complicated curves like the chicane in the lower left corner of Figure~\ref{obFree}. In this curve combination, the MPCC achieves velocities not smaller than $1$\,m/s, while the HRHC has to reduce the velocity to nearly $0.5$\,m/s to drive through the chicane, and additionally causes the car to drive an overall longer distance.

However, the MPCC has some disadvantages compared to the HRHC. In our current implementation, the MPCC still tracks the center line, which can be seen on the first straight, after the top left corner where a movement to the center line is clearly visible in Figure~\ref{obFree}. This is clearly due to the objective function of the MPCC, and cannot be completely avoided, even if the contouring cost is small. The HRHC on the other hand does not have any cost related to the center line and only uses it as a measure of progress, thus the car drives straight in the aforementioned segment of the track. The second disadvantage can be seen in the S-curve at lower end of the track, where the MPCC has a small S shape in the driven trajectory. The HRHC on the other hand shows that the section can be driven completely straight, which gives a speed benefit. This can be most probably also attributed to the contouring error penalty, for which it would be beneficial to drive an S-shape in order to follow the center line.

In summary, our results indicate that the HRHC is limited by the path planner, and the fact that the whole horizon has the same velocity. Allowing multiple velocities within one prediction horizon would most probably increase the performance of the controller. However, the complexity in the path planner grows exponentially with the number of different velocity triples in the horizon. The MPCC, which does not have this limitation, is able to plan better trajectories which also improve the closed loop performance. However, the MPCC comes at a price, first in terms of computational cost (it is about a factor of five more expensive than the HRHC when comparing computation times on the same platform, see~\secref{sec:results:compTimes} for more details on computation times) and, second, in terms of robustness. During our experiments, the re-linearization scheme turned out to be sensitive to measurement errors and model drift. In particular, sudden unexpected skidding of the car can cause problems in the re-linearization scheme. Such drifts are also responsible for the high variations between different trajectories in certain areas of the track.  Furthermore, since the MPCC relies on the results from the previous iteration to compute the linearizations, sudden jumps in the high-level corridor path planner might make the previous trajectory infeasible. In such a situation, several time steps might be needed to recover feasibility of the planned trajectory, although by using the soft constrained formulation sensible inputs are still provided to the car. The HRHC on the other hand does not need any information from the last iteration as only the current position and velocity is needed. This helps to deal efficiently with unexpected behavior and fast changing measurements or corridors from the high-level obstacle avoidance mechanism.

\subsection{Racing with multiple cars}

Since both controllers plan a path around opponents while maximizing the progress, fast and safe overtaking maneuvers are enabled. In the current work, we focus on static obstacles, which seems sufficient to outperform non-expert human drivers. To robustly overtake dynamically moving obstacles (other automatic controllers or expert drivers), it would be necessary to have a prediction of their behavior, which is currently not available. These predictions could however be systematically incorporated into the high-level corridor planner outlined in \secref{sec:dp}, and would merely change the corridor issued to the controllers presented in \secref{sec:hrhc} and \secref{sec:mpcc}.

\begin{figure}
\centering
\footnotesize
HRHC controller (\secref{sec:hrhc}) \hspace{3.75cm} MPCC controller (\secref{sec:mpcc})\\
\normalsize
\includegraphics[width = \textwidth, trim= 0cm 0cm 0cm 1.2cm, clip]{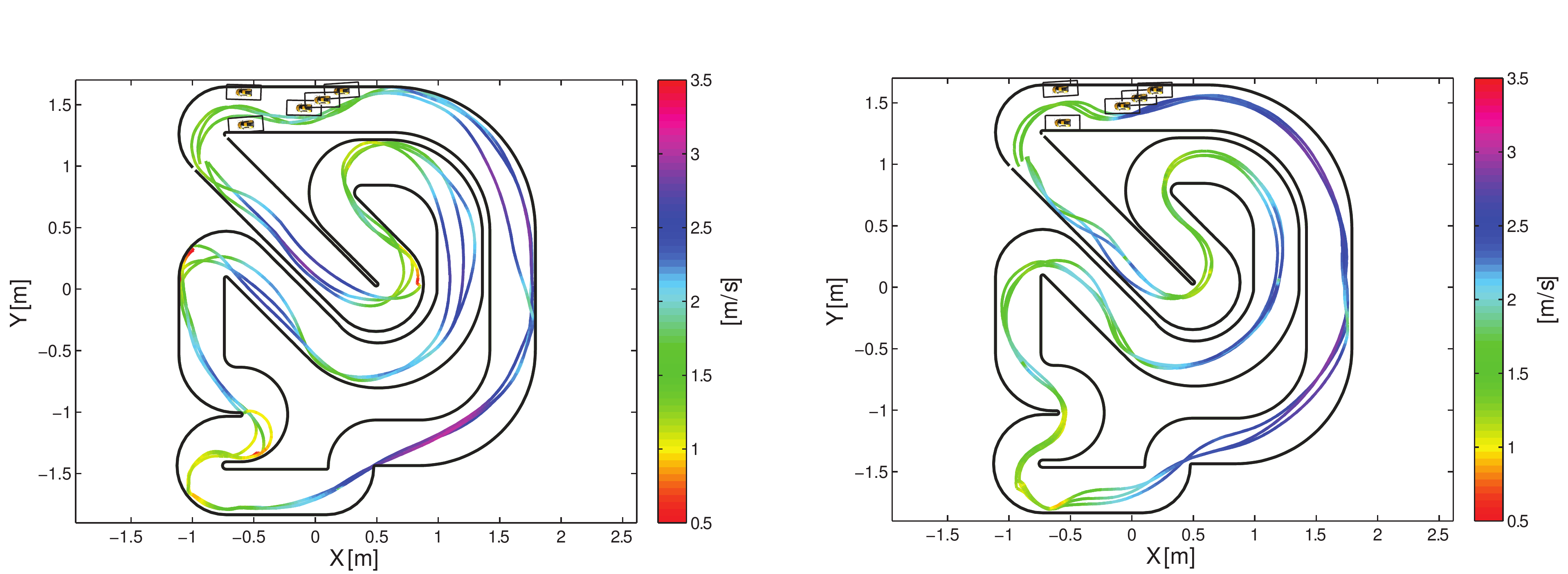}
  \caption{Driven trajectories with static obstacles for 3 different laps. A video of the MPCC with obstacle avoidance is available at \url{https://www.youtube.com/watch?v=mXaElWYQKC4}.}
  \label{obAvoid}
\end{figure}

\begin{figure}
\centering
\footnotesize
HRHC controller (\secref{sec:hrhc}) \hspace{3.75cm} MPCC controller (\secref{sec:mpcc})\\
\normalsize
\includegraphics[width = \textwidth, trim= 0cm 0cm 0cm 0cm, clip]{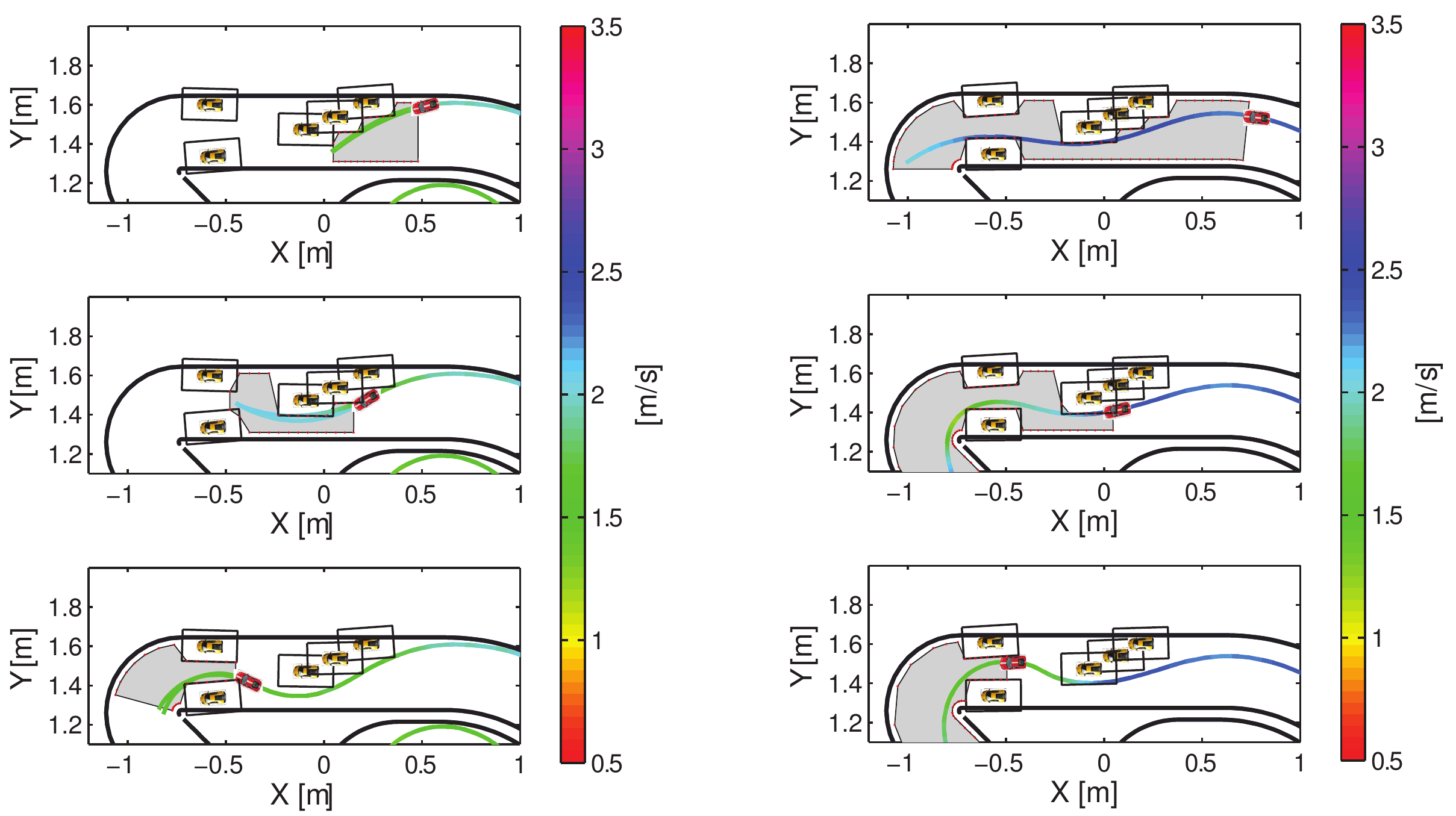}
  \caption{Obstacle avoidance at 3 different time steps. Grey indicates the corridor chosen by the obstacle avoidance algorithm.}
  \label{obAvoidZoom}
\end{figure}

An obstacle avoidance situation is shown in Figure~\ref{obAvoid}, with a detailed zoom taken at three different time points in Figure~\ref{obAvoidZoom}. These zooms show the prediction horizon and the corresponding planned trajectories. Due to the limited path planner, the HRHC is not able to directly plan a path around all obstacles from the beginning. Thus it avoids the first three cars, without predicting how to overtake the next two cars. Only after the car has made enough progress and overtaken the first group of obstacles, can the controller find a way between of the other two cars. The MPCC on the other hand plans the path around all obstacles before it even reaches the first car due to its long horizon. This is an advantage, but if the model mismatch is significant it can lead to problems, as the planned path is too optimistic.

\subsection{Computation Times} \label{sec:results:compTimes}

The computation times of the most time-consuming components of the two controllers are given in Tables~\ref{compTimeHRHC} and~\ref{compTimeMPCC} for the laps depicted in Figure~\ref{obFree} (single car racing) and Figure~\ref{obAvoid} (avoidance of opponents). Note that the computation times of the individual blocks of the two controllers cannot be compared directly. To directly compare the computation times of the two controllers, the horizon length should be identical. However, due to the constant velocity limitation in the HRHC path planner, a horizon length identical to the MPCC is impossible without using multiple segments of constant velocity, which would make the path planning combinatorial in complexity and thus currently prohibitive for a real-time implementation. Vice versa, the MPCC cannot work robustly with the short horizon length used by the HRHC. The multiple car case does not have a significant impact on the execution times of the lower level controller; only the path planning in the HRHC gets slightly more complicated. However, the computation time of the high-level corridor planner (the DP approach from \secref{sec:dp}) varies significantly, depending on the complexity of the racing situation. For the MPCC, the limiting factor with respect to the sampling time is the computation time for solving the QP. The main bottleneck in the HRHC is the path planner, which has a very high maximal computation time caused by back up rules in the case no feasible trajectory can be found. The QP on the other hand is the most expensive step in the average, but does not have a large variation in computation time, which makes it less critical. 

The $20$\,ms sampling time was missed by the HRHC controller in only $0.07$\,\% (one sampling instant in three laps) and by the MPCC controller in $4.4$\,\% (60 sampling instants in three laps) of the time. Overall, our experiments demonstrate that both schemes can be implemented in (soft) real-time at sampling rates of $50$\,Hz.

\begin{table}
\caption{HRHC computation times in milliseconds}\label{compTimeHRHC}
\begin{minipage}[t]{0.48\textwidth}
\centering
\ra{1.3}
\begin{tabular}{l r | r | r}
\toprule
Without Obstacles & Mean & Stdev & Max\\
\midrule
Border Adjustment 	& 0.02 	& 0.01 	& 0.04\\
Path Planning 		& 2.93 	& 1.96  	& 13.47\\
QP Generation 		& 0.52	& 0.09 	& 1.32\\
QP with FORCES 	& 5.10 	& 0.73 	& 7.56\\
\bottomrule
\end{tabular}
\end{minipage}
\hfill
\begin{minipage}[t]{0.48\textwidth}
\centering
\ra{1.3}
\begin{tabular}{l r | r | r}
\toprule
With Obstacles & Mean & Stdev & Max\\
\midrule
Border Adjustment 	& 1.33	& 0.37 	& 2.45 	\\
Path Planning 					& 3.18 	&  2.04  	& 9.76	\\
QP Generation 					& 0.53 	&  0.06 	& 0.76	\\
QP with FORCES 				& 5.11 	&  0.92	& 9.25	\\
\bottomrule
\end{tabular}
\end{minipage}
\end{table}

\begin{table}
\caption{MPCC computation times in milliseconds}\label{compTimeMPCC}
\begin{minipage}[t]{0.48\textwidth}
\centering
\ra{1.3}
\begin{tabular}{l r | r | r}
\toprule
Without Obstacles & Mean & Stdev & Max\\
\midrule
Border Adjustment 	& 0.34 	& 0.11 & 0.64 \\
QP Generation 		& 1.95 	& 0.44 & 2.92\\
QP with FORCES 	& 15.03 	& 1.61 & 21.16\\
\bottomrule
\end{tabular}

\end{minipage}
\hfill
\begin{minipage}[t]{0.48\textwidth}

\centering
\ra{1.3}
\begin{tabular}{l r | r | r}
\toprule
With Obstacles & Mean & Stdev & Max\\
\midrule
Border Adjustment 	& 0.70 	& 0.24 	& 1.05 \\
QP Generation 		& 2.34 	& 0.37 	& 2.79\\
QP with FORCES 	& 14.43 	& 1.40 	& 21.46\\
\bottomrule
\end{tabular}
\end{minipage}
\end{table}

\section*{Acknowledgment}

We gratefully acknowledge the work of our students in the course of the race car project: Kenneth Kuchera, Samuel Zhao and Florian Perrodin for their fruitful help in implementing the MPCC approach; Michael Janser for his contribution to the obstacle avoidance (DP) algorithm; Sandro Merkli, Nils Wenzler and Marcin Dymczyk for the development of the embedded control software and the testbed setup; Celestine D\"unner and Sandro Merkli for building the track and low-level filtering software; Michael Dahinden and Christian Stocker for the initial design of the embedded PCB board; Benjamin Keiser for the current version of the infrared vision system. We furthermore thank Sean Summers, Nikolaos Kariotoglou, Christian Conte and John Lygeros from the Automatic Control Laboratory for supervising various projects and providing important comments to the methods presented in this paper. Special thanks go to Colin Jones, who originated the idea of the race car testbed. 

The research leading to these results has received funding from the EU in FP7 via EMBOCON (ICT-248940).

\bibliographystyle{ieeetr}
\bibliography{AutonomousRacing}

\end{document}